\documentclass[11pt]{article}
\usepackage{amsmath,amssymb,epsf,wrapfig}
\usepackage[dvipdfm]{graphicx}

\newtheorem{definition}{Definition}[section]
\newtheorem{theorem}[definition]{Theorem}
\newtheorem{corollary}[definition]{Corollary}
\newtheorem{remark}[definition]{Remark}
\newtheorem{example}[definition]{Example}
\newtheorem{lemma}[definition]{Lemma}
\newtheorem{proposition}[definition]{Proposition}

\newcommand{\qed}{\hbox{\rule[-2pt]{3pt}{6pt}}}
\newenvironment{proof}{{\it Proof.}}{\hfill \qed}

\newcommand{\Z}{{\mathbb Z}}

\begin{document}
\title{Biquandles with structures related to \\ virtual links and twisted links}
\author{Naoko Kamada \\ 
Graduate School of Natural Sciences, Nagoya City University, \\ 
Aichi, 467-8501, Japan\\ 
kamada@nsc.nagoya-cu.ac.jp\\ 
\\ 
and \\ 
\\ 
Seiichi Kamada \\ 
Department of Mathematics, Hiroshima University, \\
Hiroshima 739-8526, Japan\\
kamada@math.sci.hiroshima-u.ac.jp}
\maketitle

\abstract{We introduce two kinds of  structures, called v-structures and t-structures,  on biquandles.  These structures are used for colorings of diagrams of  virtual links and twisted links such  that the numbers of colorings are invariants. 
Given a biquandle or a quandle, we give a method of constructing a biquandle with these structures.  Using the numbers of colorings, we show that Bourgoin's twofoil and non-orientable virtual $m$-foils 
do not represent virtual links.} 

\vspace{0.5cm} 

\vspace{0.5cm} 

\section{Introduction}

A {\em virtual link diagram} is an oriented link diagram possibly with encircled crossings, called virtual crossings, that are neither  positive crossings nor negative crossings.  Two diagrams are {\em equivalent} if there is a sequence of the generalized Reidemeister moves defined in \cite{Kauf99}, which are generated by moves R1, \dots, R3, V1, \dots, V4 in Figure~\ref{fgtwistmoves}.   The equivalence class of a virtual link diagram is called a {\em virtual link}.  Virtual links correspond to stable equivalence classes of oriented links in the trivial $I$-bundles over closed  orientable surfaces \cite{CKS02, KK00, Ku03}.   A {\em twisted link diagram} is a virtual link diagram which may have some bars on edges.   Two diagrams are {\em equivalent} if there is a sequence of the extended Reidemeister moves defined in \cite{Bo08}, 
which are generated by all moves  in Figure~\ref{fgtwistmoves}.   The equivalence class of a twisted link diagram is called a {\em twisted link}.   
(In \cite{Bo08} the extended Reidemeister moves are illustrated without orientations.  Note that all moves with possible orientations are obtained by combining the moves in Figure~\ref{fgtwistmoves}.  For example, see Figure~\ref{fgtwistmovesB}.)   
Twisted links correspond to stable equivalence classes of oriented links in oriented $3$-manifolds that are orientation $I$-bundles over closed but not necessarily orientable surfaces \cite{Bo08}.  

In this paper we define two kinds of structures on biquandles which are related to virtual links and twisted links.   

{\begin{figure}[h]\begin{center}
\includegraphics[width=10.0cm,clip]{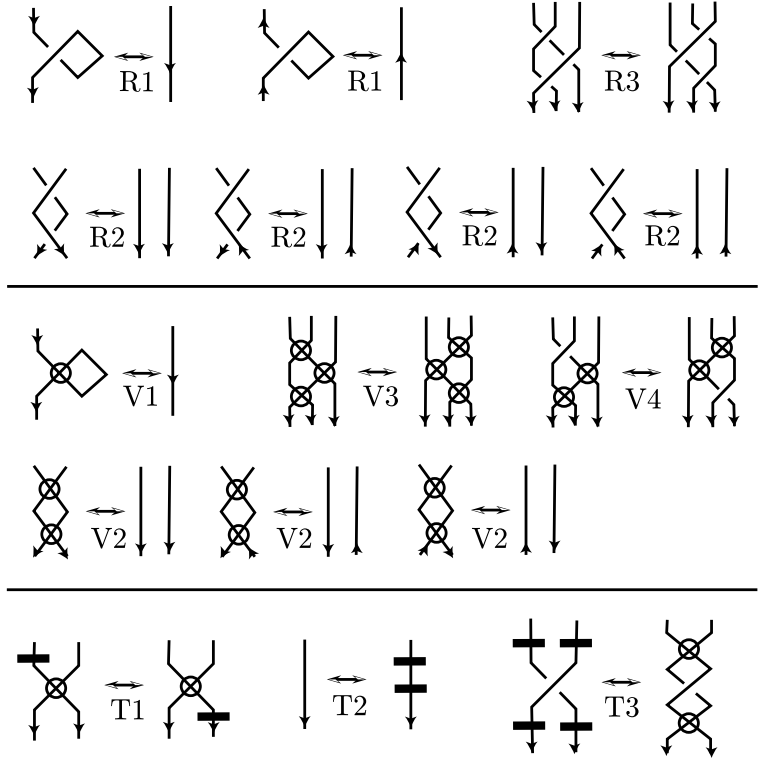}\vspace{0cm}
\caption{The extended Reidemeister moves}\label{fgtwistmoves}
\end{center}\end{figure}}

{\begin{figure}[h]\begin{center}
\includegraphics[width=6.0cm,clip]{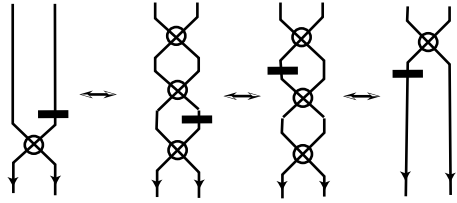}\vspace{0cm}
\caption{Another T1 move}\label{fgtwistmovesB}
\end{center}\end{figure}}

A  {\em biquandle} is a pair $(X,R)$ consisting of  a set $X$ and a bijection $R: X^2 \to X^2$ satisfying certain conditions corresponding to Reidemeister moves for classical link diagrams \cite{FJK04, FRS93, St06}  (Section~\ref{sect:basicsbiquandles}).  
In Section~\ref{sect:vtbiquandles} we introduce the notions of  a v-structure $V: X^2 \to X^2$ and a t-structure $T: X \to X$ 
which are additional structures on a biquandle $(X,R)$ related to virtual links and twisted links.   The pair $(V, T)$ is called a vt-structure of $(X,R)$.  

A coloring of a virtual link diagram by a v-structured biquandle $(X, R, V)$ or  
a coloring of a twisted link diagram by a vt-structured biquandle $(X, R, V, T)$ is defined as follows:  Let $D$ be a diagram of a virtual link or a twisted link.  The {\em edges} of $D$ mean the connected arcs obtained when all the real crossings, virtual crossings and bars are removed.  

\begin{definition}\label{def:coloring}{\rm 
A {\em coloring} of $D$ by $(X, R, V)$ or  $(X, R, V, T)$ 
is a map from the set of edges of $D$ to $X$ such  that for each crossing or bar, say $v$, of $D$, if $x_1, x_2, x_3, x_4$ are elements of $X$ assigned the edges around $v$ as in Figure~\ref{fglabelreal} then 
\begin{itemize}
\item[(1)] $R(x_1, x_2) = (x_3, x_4)$ when $v$ is a positive crossing,   
\item[(2)] $R^{-1}(x_1, x_2) = (x_3, x_4)$ when $v$ is a negative crossing,  
\item[(3)] $V(x_1, x_2) = (x_3, x_4)$ when $v$ is a virtual crossing, and   
\item[(4)] $T(x_1) = x_2$ when $v$ is a bar. 
\end{itemize}  
We also call a coloring by $(X, R, V)$ an {\em $(X, R, V)$-coloring}, and a coloring by $(X, R, V, T)$ an {\em $(X, R, V, T)$-coloring}. 
}\end{definition}

The concept of a coloring by $(X, R, V)$ in Definition~\ref{def:coloring} and  the following theorem (Theorem~\ref{thm:coloringvirtual}) were considered in \cite{BaF11}.  
Refer to \cite{BaF11} for examples.  

{\begin{figure}[h]\begin{center}
\includegraphics[width=8.0cm,clip]{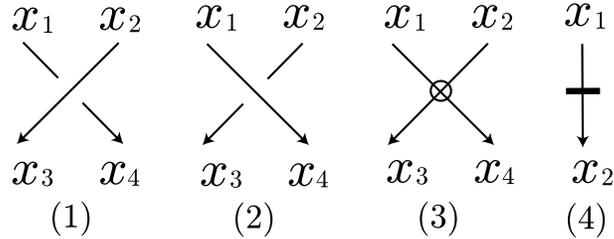}\vspace{0cm}
\caption{Colorings}\label{fglabelreal}
\end{center}\end{figure}}

\begin{theorem}[\cite{BaF11}]\label{thm:coloringvirtual}
If $D$ and $D'$ are  virtual link diagrams representing the same virtual link, then there is a bijection between the set of colorings of $D$ by a v-structured biquandle $(X, R, V)$ and that of $D'$.   
\end{theorem} 

This is generalized to twisted link diagrams. 

\begin{theorem}\label{thm:coloring}
If $D$ and $D'$ are  twisted link diagrams representing the same twisted link, then there is a bijection between the set of colorings of $D$ by a vt-structured biquandle $(X, R, V, T)$ and that of $D'$.   
\end{theorem} 

Therefore  the number of colorings by a vt-structured biquandle $(X, R, V, T)$ is an invariant of a twisted link.  

Given a biquandle $(X_0, R_0)$ and   
automorphisms  $f$
and $g$ with $f^2=1$ and $fg = gf$, 
we give a method of constructing a vt-structured biquandle $(X, R, V_f, T_g)$, which we call a 
{\em twisted product} of $(X_0, R_0)$.  

Let $(X_0, R_0)$ be a biquandle.  
We use the notation due to \cite{FJK04} such that for $a, b \in X_0$, 
$$R_0 (a,b) = (b_a, a^b),$$ 
namely, $b_a = p_1 R_0(a,b)$ and $a^b = p_2 R_0(a,b)$, where $p_i : X_0 \times X_0 \to X_0$ is the $i$th factor projection.    

\begin{theorem}\label{thm:twistedproduct}
Let $(X_0, R_0)$ be a biquandle.  
Let $X = X_0 \times X_0$.  Define a map $R : X^2 \to X^2$ by 
$$ R( (a_1, b_1), (a_2, b_2)) = (( {a_2}_{a_1}, {b_2}^{b_1}),  ({a_1}^{a_2}, {b_1}_{b_2})). $$ 
For automorphisms $f$ and $g$ of $(X_0, R_0)$, 
define maps $V_f: X^2 \to X^2$ and $T_g: X \to X$  by 
\begin{eqnarray*}
V_f ( (a_1, b_1), (a_2, b_2) ) &=& ( (f^{-1}a_2, f^{-1} b_2), ( f a_1, f b_1 )),  \mbox {and} \\  
T_g(a, b) &=& (g^{-1} b, g a). 
\end{eqnarray*}
Then the following holds.  

$(1)$ $(X, R)$ is a biquandle.  

$(2)$ $V_f$ is a  v-structure of $(X, R)$.  

$(3)$ Suppose that $f^2 =1$ and $fg =gf$.  Then $(V_f, T_g)$ is a vt-structure of $(X, R)$. 
\end{theorem} 

\begin{definition} {\rm 
In the situation of Theorem~\ref{thm:twistedproduct}, we call the biquandle  $(X, R)$ the  {\em twisted product biquandle} of $(X_0, R_0)$, and  the vt-structured biquandle 
$(X, R, V_f, T_g)$ a {\em twisted product} of $(X_0, R_0)$.   When $f=g=1$, we call the quadruplet  the {\em standard twisted product} of $(X_0, R_0)$.  
}\end{definition} 

A {\em quandle} is a pair $(Q, \ast)$ consisting of a set $Q$ and a binary operation $\ast: Q \times Q \to Q$, $(a, b) \mapsto a \ast b$, such that  (i) for any $a \in Q$, $a \ast a=a$, (ii) for any $a, b \in Q$, there exists a unique element $c$ with $c \ast b =a$, and (iii) for any $a, b , c \in Q$, $(a \ast b) \ast c = (a \ast c) \ast (b \ast c)$ \cite{FR92, Joyce1982, Matveev82}.   The {\em dual} operation $\overline{\ast}$ of $\ast$ is a binary operation $\overline{\ast}: Q \times Q \to Q$ such that $a \,\overline{\ast}\, b = c \Longleftrightarrow c \ast b =a $.   In Section~\ref{sect:proofs}, we use Fenn and Rourke's notation   \cite{FR92}: $a \ast b$ and  $a \,\overline{\ast}\, b$ are denoted by $a^b$ and  $a^{\overline b}$ (or $a^{b^{-1}}$), respectively, and 
$a^{bc}$ means $(a^b)^c$, etc.   

When $(X_0, R_0)$ is the  {\em biquandle derived from a quandle $(Q, \ast)$}, i.e.,  
$X_0=Q$ and  $ R_0 (x, y) = (y, x \ast y)$, we call the twisted product biquandle $(X, R)$ of $(X_0, R_0)$ the {\em twisted product  biquandle of $(Q, \ast)$}. 
If $f$ and $g$ are quandle automorphisms of $(Q, \ast)$, then they are biquandle automorphisms of $(X_0, R_0)$.  Suppose that $f^2=1$ and $fg = gf$.  
A {\em twisted product of $(Q, \ast)$} means a twisted product $(X, R, V_f, T_g)$ of $(X_0, R_0)$.  
The operations $R, V_f, T_g$ on $X=Q^2$ are given as follows: 
\begin{eqnarray*}
R( (a_1, b_1), (a_2, b_2) ) &=& (( a_2, {b_2} \ast {b_1}), ( {a_1} \ast {a_2}, b_1)), \\ 
V_f ( (a_1, b_1), (a_2, b_2) ) &=& ( (f a_2, f b_2), ( f a_1, f b_1 )),  \mbox { and} \\ 
T_g(a, b) &=& (g^{-1} b, g a). 
\end{eqnarray*}
When $f=g=1$, we have the following definition.  

\begin{definition}\label{def:twistedproductquandleA} {\rm 
The {\em standard twisted product} of a quandle $(Q, \ast)$, which we denote by ${\cal B}(Q, \ast)$,  is a vt-structured biquandle $(X, R, V, T)$ such that  
$X= Q^2$ and 
\begin{eqnarray*}
R( (a_1, b_1), (a_2, b_2) ) &=& (( a_2, {b_2} \ast {b_1}), ( {a_1} \ast {a_2}, b_1)), \\ 
V ( (a_1, b_1), (a_2, b_2) ) &=& ( (a_2, b_2), ( a_1, b_1 )),  \mbox {and} \\  
T(a, b) &=& (b, a). 
\end{eqnarray*}
}\end{definition}

For a quandle $(Q, \ast)$, the number of upper $(Q, \ast)$-colorings 
of a virtual link diagram is an invariant of a virtual link, and so is  
that of lower $(Q, \ast)$-colorings \cite{Kauf99}. (A geometric interpretation of the upper/lower knot quandles of a virtual link is given in \cite{KK00}.  The upper/lower $(Q, \ast)$-colorings correspond to the homomorphisms from these geometric quandles to $(Q, \ast)$ as in the classical case \cite{FR92, Joyce1982, Matveev82}.) 

Let ${\cal B}(Q, \ast)$ be the standard  twisted product of a quandle $(Q, \ast)$.  

\begin{theorem} \label{thm:bourgoinproductcoloring} 
If $D$ is a twisted link diagram which is equivalent to a virtual link  diagram $D'$, then 
the number of ${\cal B}(Q, \ast)$-colorings of $D$ is the product of the number of upper $(Q, \ast)$-colorings 
of $D'$ and that of lower $(Q, \ast)$-colorings 
of $D'$.  
\end{theorem} 

\begin{corollary}\label{cor:bourgoinproductcoloring}
Let $(Q, \ast)$ be a finite quandle with $n$ elements.  
If the number of ${\cal B}(Q, \ast)$-colorings of a twisted link diagram $D$  is less  than $n^2$, then $D$ does not represent a virtual link.    
\end{corollary} 

\proof \quad Every virtual link diagram has at least 
$n$ trivial upper $(Q, \ast)$-colorings  and at least $n$ trivial lower $(Q, \ast)$-colorings.  If $D$ is equivalent to a virtual link diagram,  by Theorem~\ref{thm:bourgoinproductcoloring} there are at least $n^2$ colorings of $D$ by ${\cal B}(Q, \ast)$.  \qed  

\vspace{0.3cm}
Using an argument due to \cite{Do90}, Bourgoin \cite{Bo08} showed that the twisted link diagram illustrated in Figure~\ref{fgtwistedtrefoil}, which we call {\em Bourgoin's twofoil},  does not represent a classical link.   He also defined a twisted link invariant, called the {\em twisted Jones polynomial},  and 
showed that Bourgoin's twofoil has the same twisted Jones polynomial with a certain virtual link diagram (Figures 6 and 7 of \cite{Bo08}).  Thus, as mentioned in \cite{Bo08},  one cannot distinguish it from virtual links by use of the twisted Jones polynomial.  

{\begin{figure}[h]\begin{center}
\includegraphics[width=4.0cm,clip]{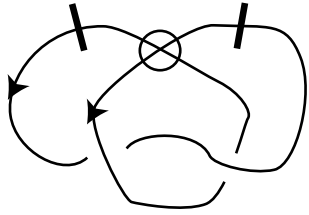}\vspace{0cm}
\caption{Bourgoin's twofoil}\label{fgtwistedtrefoil}
\end{center}\end{figure}}

We call a diagram on the left-hand side of Figure~\ref{fgtwistedmfoil} {\em a non-orientable virtual $m$-foil} and denote it by $F_m$, 
where $m$ is the number of the real crossings $(m \geq 1)$.  When $m=2$, it is Bourgoin's twofoil.  
It has a realization as a link diagram on a projective plane depicted 
on the right-hand side of the figure, where the unit disk is made into a projective plane.   
By a calculation using induction on $m$, we see that the twisted Jones polynomial of $F_m$ is 
$$ A^{-2m} ( A^{-4m} + (-1)^{m+1} ( 1 + A^2 + A^{-2})) \in {\bf Z}[A^{\pm 1}, M].  $$
This polynomial is also the twisted Jones polynomial of a virtual link diagram obtained when the two bars are removed from $F_m$.  Thus one cannot distinguish it from virtual links by use of the twisted Jones polynomial. 

{\begin{figure}[h]\begin{center}
\includegraphics[width=8cm,clip]{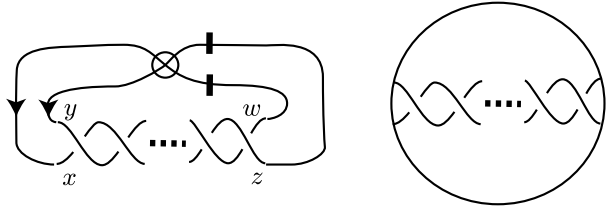}\vspace{0cm}
\caption{non-orientable virtual $m$-foil}\label{fgtwistedmfoil}
\end{center}\end{figure}}

In Section~\ref{sect:proofs} we study colorings of the diagram $F_m$ by the standard twisted product ${\cal B}(Q, \ast)$, and show the following.     

\begin{theorem}\label{thm:twofoil}
\begin{itemize}
\item[$(1)$] For $m > m'  \geq 1$, $F_m$ and $F_{m'}$ represent distinct twisted links.  
\item[$(2)$] 
For $m \geq 1$, $F_m$ does not represent a virtual link.  
\end{itemize}
\end{theorem} 

The first assertion of this theorem is also seen by the twisted Jones polynomials. 
Now we have an infinite family of twisted links which are not virtual links,  but they are not distinguished from virtual links by the twisted Jones polynomials.  This example was suggested the authors by Roger Fenn. 

\begin{remark}\label{remark:twistedproductquandleB} {\rm 
Given a quandle $(Q, \ast)$, let  
$X= Q^2$ and 
\begin{eqnarray*}
R( (a_1, b_1), (a_2, b_2) ) &=& (( a_2, {b_2} \, \overline{\ast}\,  {b_1}), ( {a_1} \ast {a_2}, b_1)), \\ 
V ( (a_1, b_1), (a_2, b_2) ) &=& ( (a_2, b_2), ( a_1, b_1 )),  \mbox {and} \\  
T(a, b) &=& (b, a). 
\end{eqnarray*}
Then $(X, R)$ is a biquandle, and $V$ is a v-structure of $(X, R)$.  Since 
\begin{eqnarray*}
(T \times T)R(T\times T) ( (a_1, b_1), (a_2, b_2) ) &=& ( (a_2  \, \overline{\ast}\, a_1, b_2), (a_1, b_1 \ast b_2) ) \quad \mbox{ and } \\ 
 VRV ( (a_1, b_1), (a_2, b_2) ) &=& ( (a_2 \ast a_1, b_2), (a_1, b_1 \, \overline{\ast}\, b_2) ), 
\end{eqnarray*}
the operation $T: X^2 \to X^2$ does not satisfy $(T \times T)R(T\times T) = VRV$ unless $(Q, \ast)$ is an {\em involutory}  quandle, i.e., $\ast = \overline \ast$.  
Thus the quadruplet 
$(X, R, V, T)$ is not a vt-structured biquandle and one should not use this for colorings  of  twisted links.  
}\end{remark} 

We recall the notion of a biquandle, and prepare some lemmas in 
Section~\ref{sect:basicsbiquandles}.   In Section~\ref{sect:vtbiquandles} the definitions of a t-structure and a v-structure are given, and Theorem~\ref{thm:twistedproduct} is proved.  
Section~\ref{sect:proofs} is devoted to proofs of  Theorems~\ref{thm:coloring}, \ref{thm:bourgoinproductcoloring} and \ref{thm:twofoil}.

\section{Biquandles}\label{sect:basicsbiquandles}

For a set $X$ we denote by $X^n$ the $n$-fold Cartesian product of  $X$, and denote by $p_i: X^n \to X$ the $i$th factor projection  for each $i =1, \dots, n$.   The composition $g \circ f$ of two maps $f$ and $g$ is also denoted by $g \cdot f$ or $gf$.   The identity map $x \mapsto x$ on $X$ is denoted by $1_X$ or $1$, and 
the transposition map $(x, y) \mapsto (y,x)$ on $X^2$ is denoted by $\tau_{X^2}$ or $\tau$.  

The basic idea of a birack was given in \cite{FRS93}.  
The following is the definition of a (strong) birack and a (strong) biquandle introduced by R.~Fenn, M.~Jordan-Santana and L.~Kauffman (Definitions~4.2 and 4.6 of \cite{FJK04}).  

\begin{definition}\label{def:biquandlefenn}{\rm (\cite{FJK04}) 
A pair $(X, R)$ of a set $X$ and a bijection $R: X^2 \to X^2$ is a {\em birack} if the following conditions 
(B1) and (B2) 
are satisfied.   It is a  {\em biquandle} if  (B1), (B2) and (B3) are    satisfied. 
\begin{itemize}
\item[(B1)] $R$ satisfies the set-theoretic Yang-Baxter equation, i.e., 
$$ (R \times 1) (1 \times R) (R \times 1) = (1 \times R)  (R \times 1)  (1 \times R) : X^3 \to X^3.$$ 
\item[(B2)] For  $a, b \in X$, let $f_a: X \to X$ and $f^b: X \to X$ be maps defined by 
$$f_a(x) = p_1 R(a,x) \quad \mbox{and} \quad f^b(x) = p_2 R(x,b).$$
Then both $f_a$ and $f^b$ are bijections for every $a$ and $ b$ of $X$.  
\item[(B3)] For every $a$ and $b$ of $X$, 
$$    (f_a)^{-1}(a) = f^{(f_a)^{-1}(a)} (a)  \quad   \mbox{and}  \quad (f^b)^{-1}(b) = f_{(f^b)^{-1}(b)} (b).  $$
\end{itemize}
}\end{definition}

As in Section 3 of  \cite{FJK04} a birack/biquandle operation $R: X^2 \to X^2$  defines two binary operations on $X$; $(a,b) \mapsto b_a$ and $(a,b) \mapsto a^b$  such that $R(a,b) =(b_a, a^b)$.   Refer to  \cite{BaF08, BuF04, CSWEM09, F09, FJK04, HK07, KauMan05} for examples of biracks and biquandles.  

Let $(X, R)$ and $(X', R')$ be biracks or biquandles.  A map $h: X \to X'$ is called a {\em homomorphism} if 
$( h \times h)   R = R'  (h \times h) : X^2 \to X^2$.   We denote it by $h: (X, R) \to (X', R')$.   An {\em isomorphism} is a bijection which is a homomorphism. 

Using the notion of a sideways operation due to Fenn, et al. \cite{BuF04, F09, FJK04}, we can restate Definition~\ref{def:biquandlefenn} as follows.  

\begin{definition}\label{def:biquandle} {\rm 
A pair $(X, R)$ consisting of a set $X$ and a bijection $R: X^2 \to X^2$ is a {\em birack} if the following conditions 
(B1) and (B2$'$) 
are satisfied.   It is a  {\em biquandle} if  (B1), (B2$'$) and (B3$'$) are    satisfied. 
\begin{itemize}
\item[(B1)] 
$ (R \times 1) (1 \times R) (R \times 1) = (1 \times R)  (R \times 1)  (1 \times R) : X^3 \to X^3$. 
\item[(B2$'$)] There is a unique bijection $S : X^2 \to X^2$ such that for any $x_1, \dots,  x_4 \in X$, 
$$ S(x_1, x_3) = (x_2, x_4) \Longleftrightarrow R(x_1, x_2) = (x_3, x_4). $$
\item[(B3$'$)] There is a bijection $s: X \to X$ such that for any $x \in X$, 
$$ R(x, s(x) ) = (x, s(x) ).  $$
\end{itemize}
}\end{definition}

We call the bijections $S: X^2 \to X^2$ and $s: X\to X$  above the {\em sideways operation}  of $R$ and the {\em shift operation} of $R$, 
and denote them by ${\rm side} R$ and ${\rm shift} R$, respectively.  (Note that if  $({\rm B}2')$ holds then 
a bijection $s: X \to X$  in $({\rm B}3')$ is unique if there exists, since $S$ is unique.)  Refer to \cite{BuF04, F09, FJK04} for sideways operations.

First we observe that Definitions~\ref{def:biquandlefenn} and \ref{def:biquandle} are equivalent, and give some lemmas on biquandles. 

\begin{theorem} 
$(1)$  
For a bijection $R: X^2 \to X^2$, the conditions $({\rm B}2)$ and $({\rm B}2')$ are equivalent.  
$(2)$  
For a bijection $R: X^2 \to X^2$ satisfying $({\rm B}2)$,  the conditions $({\rm B}3)$ and $({\rm B}3')$ are equivalent.  
\end{theorem}

\proof  \quad 
(1) 
Suppose $({\rm B}2)$.  
Define maps $S: X^2 \to X^2$ and $S^{-1}: X^2 \to X^2$  by 
\begin{eqnarray*}
S(x, y ) &=& ( (f_x)^{-1}(y),   p_2 R(x, (f_x)^{-1}(y)) )  \quad \mbox{and} \quad  \\ 
S^{-1}(x,y) &=& ( (f^x)^{-1}(y),  p_1 R((f^x)^{-1}(y), x)    ).
\end{eqnarray*}
Then $S S^{-1} = S^{-1} S= 1$, and $ S(x_1, x_3) = (x_2, x_4) \Longleftrightarrow R(x_1, x_2) = (x_3, x_4)$. 
Let  $S': X^2 \to X^2$ be another bijection such that $ S'(x_1, x_3) = (x_2, x_4) \Longleftrightarrow R(x_1, x_2) = (x_3, x_4)$.  For any $x, y \in X$, since $f_{x}$ is bijective,  we have $p_1 S(x,y) = p_1 S'(x,y)$.  Then 
$p_2 S(x,y) =  p_2 R(x, p_1 S(x,y))  = p_2 R(x, p_1 S'(x,y))   = p_2 S'(x,y)  $.  Thus $S=S'$. 

Suppose $({\rm B}2')$.  
The inverse maps of $f_a$ and $f^b$ are obtained by 
$$(f_a)^{-1}(x) = p_1 S(x,a) \quad  \mbox{and} \quad (f^b)^{-1}(x) = p_1 S^{-1}(b,x).$$

(2) 
Suppose $({\rm B}3)$.  
Let $s: X \to X$ and $s^{-1}: X \to X$ be maps defined by 
$$s(x) = (f_x)^{-1}(x) \quad \mbox{and} \quad s^{-1}(y) = (f^y)^{-1}(y). $$  
Since $s(x) = (f_x)^{-1}(x) = f^{(f_x)^{-1}(x)}(x) = 
f^{s(x)}(x)$,  we have $R(x, s(x)) = (x, s(x))$.  Since 
$s^{-1}(y) =  (f^y)^{-1}(y) = f_{(f^y)^{-1}(y)} (y) = f_{s^{-1}(y)}(y) $,  we have 
$R(s^{-1}(y), y) = (s^{-1}(y), y)$.   Then $s: X \to X$ and $s^{-1}: X \to X$ are the inverse maps of each other. 

Suppose $({\rm B}3')$.  
Since $R(a, s(a)) = (a, s(a))$, we have 
$   (f_a)^{-1}(a) = s(a) = f^{(f_a)^{-1}(a)} (a)$.   Since $R(s^{-1}(b) , b ) =(s^{-1}(b) ,b )$, we have 
$ (f^b)^{-1}(b) = s^{-1}(b) = f_{(f^b)^{-1}(b)} (b)$.   
 \qed

Now one may use Definition~\ref{def:biquandle} instead of Definition~\ref{def:biquandlefenn}.  

We give some lemmas on biquandles, which are also valid for biracks.

\begin{lemma}
$(1)$ If $(X, R)$ is a biquandle, then $(X, R^{-1})$ is a biquandle with  
${\rm side} R^{-1} = \tau \cdot {\rm side} R \cdot \tau$ and 
${\rm shift} R^{-1} = {\rm shift} R $.  
\end{lemma}

\proof \quad 
Since $R^{-1} \times 1 = ( R \times 1)^{-1} : X^3 \to X^3$,  the map $R^{-1}$ satisfies the set-theoretic Yang-Baxter equation. 
Put $S' = \tau \cdot {\rm side} R \cdot \tau$.  Then 
$S'(x_3, x_1) = (x_4, x_2)$ $\Longleftrightarrow$ ${\rm side} R(x_1, x_3) = (x_2, x_4)$ $\Longleftrightarrow$ $R(x_1, x_2) = (x_3, x_4)$ 
$\Longleftrightarrow$ $(x_1, x_2) = R^{-1}(x_3, x_4)$.  The uniqueness of $S'$ follows from that of ${\rm side} R$.  
Let $s={\rm shift} R$.   Since $R(x, s(x)) = (x, s(x))$, we have $R^{-1}(x, s(x)) = (x, s(x))$.  
\qed

\begin{lemma}\label{lem:tauRtau}
 If $(X, R)$ is a biquandle, then $(X, \tau R \tau )$ is a biquandle  with  
${\rm side} (\tau R \tau) = ({\rm side} R)^{-1}$  and 
${\rm shift} (\tau R \tau) = ({\rm shift} R)^{-1} $.  
\end{lemma}

\proof \quad 
It is obvious that $ \tau R \tau$ is bijective.  
For simplicity we denote by $\tau_1, \tau_2, R_1$ and $R_2$ the maps  
$\tau \times 1, 
1 \times \tau, R\times 1$ and $1 \times R$, respectively.   Noting that 
$\tau_1 \tau_2 \tau_1 = \tau_2 \tau_1 \tau_2$ 
and $R_1 \tau_2 \tau_1 = \tau_2 \tau_1 R_2$, we see that 
$(\tau_1 R_1 \tau_1) (\tau_2 R_2 \tau_2) (\tau_1 R_1 \tau_1) = 
\tau_1 \tau_2 \tau_1 R_1 R_2 R_1 \tau_1 \tau_2 \tau_1$ and 
$(\tau_2 R_2 \tau_2) (\tau_1 R_1 \tau_1) (\tau_2 R_2 \tau_2) = 
\tau_2 \tau_1 \tau_2 R_2 R_1 R_2 \tau_2 \tau_1 \tau_2$.  
Thus $R_1 R_2 R_1 = R_2 R_1 R_2$ implies that 
$$(\tau_1 R_1 \tau_1) (\tau_2 R_2 \tau_2) (\tau_1 R_1 \tau_1) = (\tau_2 R_2 \tau_2) (\tau_1 R_1 \tau_1)  (\tau_2 R_2 \tau_2). $$
So $\tau R \tau$ satisfies the set-theoretic Yang-Baxter equation.  
$\tau R \tau (x_1, x_2) = (x_3, x_4)$ $\Longleftrightarrow$ 
$R(x_2, x_1)= (x_4, x_3)$ $\Longleftrightarrow$ 
${\rm side} R (x_2, x_4) = (x_1, x_3)$ $\Longleftrightarrow$ 
$({\rm side} R)^{-1} (x_1, x_3) = (x_2, x_4)$.  Thus ${\rm side}(\tau R \tau) = ({\rm side} R)^{-1}$.  (The uniqueness of ${\rm side}(\tau R \tau)$  follows from that of ${\rm side} R$.) 
Let $s={\rm shift} R$.   Since $R(x, s(x))= (x, s(x))$ for every $x \in X$, we have $\tau T \tau (s(x), x) = (s(x), x)$.  Thus ${\rm shift} (\tau R \tau) = ({\rm shift} R)^{-1} $.  \qed

\begin{lemma}\label{tauhomo}
Let $(X, R)$ and $(X', R')$ be biquandles, and let $f: X \to X'$ be a map.  The following three conditions are mutually equivalent. 
\begin{itemize}
\item[${\rm (i)}$] $f$ is a homomorphism from $(X, R)$ to $(X', R')$.  
\item[${\rm (ii)}$] $f$ is a homomorphism from $(X, R^{-1})$ to $(X', R'^{-1})$.  
\item[${\rm (iii)}$]  $f$ is a homomorphism from $(X, \tau R \tau)$ to $(X', \tau R' \tau)$.  
\end{itemize}
\end{lemma}

\proof \quad 
Since $( f \times f)  R = R'  (f \times f)$ $\Longleftrightarrow$ $R'^{-1}  ( f \times f)  = (f \times f)   R^{-1}$, we have (i) $\Longleftrightarrow$ (ii).   
Since $( f \times f)  R = R'  (f \times f)$ $\Longleftrightarrow$ $  \tau ( f \times f) R \tau  =   \tau (f \times f)   R' \tau $ $\Longleftrightarrow$ $ ( f \times f) \tau R \tau  = (f \times f)  \tau R' \tau $, we have (i) $\Longleftrightarrow$ (iii).   
\qed

\begin{lemma}\label{lem:prod} 
Let $(X_1, R_1)$ and $(X_2, R_2)$ be biquandles.  
Let $R: (X_1 \times X_2)^2 \to (X_1 \times X_2)^2$ be a map defined by 
$$R( (a_1, b_1), (a_2, b_2) ) = ( (p_1 R_1(a_1, a_2),  p_1 R_2(b_1, b_2)),  (p_2 R_1(a_1, a_2), p_2 R_2(b_1, b_2)).$$
Then $(X_1 \times X_2, R)$ is a  biquandle.   Moreover, if $f_i: (X_i, R_i) \to (X_i, R_i)$ $(i=1,2)$ are homomorphisms, then 
$f_1 \times f_2: (X_1 \times X_2, R) \to (X_1 \times X_2, R)$ is a homomorphism. 
\end{lemma} 

\proof \quad 
The operation $R$ satisfies  the set-theoretic Yang-Baxter equation, since $R_1$ and $R_2$ do.  
The sideways operation $S$ of $R$ is given by 
$$S( (a_1, b_1), (a_3, b_3) ) = ( (p_1 S_1(a_1, a_3),  p_1 S_2(b_1, b_3)),  (p_2 S_1(a_1, a_3), p_2 S_2(b_1, b_3)), $$
where $S_1$ and $S_2$ are the sideways opetations of $R_1$ and $R_2$.  
The shift operation $s$ of $R$ is given by 
$$s( (a, b) ) = (s_1(a),  s_2(b)), $$
where $s_1$ and $s_2$ are the shift opetations of $R_1$ and $R_2$.  
It is obvious that $f_1 \times f_2$ is a homomorphism. 
\qed 

We call the biquandle $(X_1 \times X_2, R)$ in Lemma~\ref{lem:prod} 
the {\em direct product} of $(X_1, R_1)$ and $(X_2, R_2)$.

\section{v- and t-Strucures on biquandles}\label{sect:vtbiquandles}

First we introduce the notion of a v-structure.

\begin{definition}\label{def:tstructure}{\rm (cf. \cite{BaF11})
Let $(X, R)$ be  a biquandle. 
A bijection $V: X^2 \to X^2$ is a {\em v-structure} of $(X, R)$ if the following conditions are satisfied. 
\begin{itemize}
\item[(1)] $(X, V)$ is a biquandle. 
\item[(2)] $V^2 =1: X^2 \to X^2$. 
\item[(3)] $(V \times 1) (1 \times V) (R \times 1) = (1 \times R) (V \times 1) (1 \times V) : X^3 \to X^3$. 
\end{itemize} 
We call $(X, R, V)$ a {\em v-structured biquandle}.  
}\end{definition} 

A v-structure is used for colorings of  virtual link diagrams.   

\begin{example} {\rm 
Let $(X, R)$ be a  biquandle.  Let $V$ be the transposition $\tau: X^2 \to X^2, (x,y) \mapsto (y,x)$.  It is a v-structure.  
A biquandle with this v-structure is used for colorings of  virtual link diagrams in  \cite{CSWEM09, FJK04, HK07, KR03}, etc.   
}\end{example}

\begin{example}\label{example:virtualquandle} {\rm 
Let $(X, R)$ be  a biquandle.  Let $f: (X, R) \to (X, R)$ be an automorphism.   
Let $V: X^2 \to X^2$ be a map defined by $V(x_1,x_2) = (f^{-1} x_2, f x_1)$.  It is a v-structure.  
A  biquandle with this structure is called  a {\em virtual biquandle}.  
See \cite{KauMan05} and Definition~3.3 of \cite{CN09}.   
}\end{example}

We introduce the notion of a t-structure, or a 
vt-structure, which is related to twisted links.  

\begin{definition}\label{def:vtstructure}{\rm 
Let $(X, R, V)$ be a v-structured biquandle. 
A bijection $T: X \to X$ is a {\em t-structure} of $(X, R, V)$ if the following conditions are satisfied. 
\begin{itemize}
\item[(1)] $T^2 =1$. 
\item[(2)] $V (T \times 1) = ( 1 \times T) V$. 
\item[(3)] $(T \times T) R (T \times T) = V R V$. 
\end{itemize} 
We call $(X, R, V, T)$ a {\em vt-structured biquandle}, and   $(V, T)$ a  {\em vt-structure} of $(X, R)$.  
}\end{definition} 

Since $V^2=1$, the condition (2) of Definition~\ref{def:vtstructure} is equivalent to that $V (1 \times T) = (T \times 1) V$.  

\vspace{0.3cm}
\indent 
{\it Proof of Theorem~\ref{thm:twistedproduct}}. 
(1)  By Lemma~\ref{lem:tauRtau}, when $(X_0, R_0)$ is a biquandle, so is $(X_0, \tau R_0 \tau)$.  By Lemma~\ref{lem:prod}, 
$(X, R)$ is a biquandle, since it is the direct product of $(X_0, R_0)$ and $(X_0, \tau R_0 \tau)$.   (2) It follows from Lemmas~\ref{tauhomo} and \ref{lem:prod} that $f \times f:  X =X_0 \times X_0 \to X =X_0 \times X_0$ is an automorphism of $(X, R)$.  By Example~\ref{example:virtualquandle} we have (2).   
The assertion (3) is verified by direct calculation, which is left to the reader.   \qed

\section{Proof of Theorems} \label{sect:proofs} 

In this section we prove Theorems~\ref{thm:coloring}, \ref{thm:bourgoinproductcoloring} and \ref{thm:twofoil}. 

\vspace{0.3cm}
\indent 
{\it Proof of Theorem~\ref{thm:coloring}}. Let $\Delta$ be a $2$-disk where a move in Figure~\ref{fgtwistmoves} transforms $D$ to $D'$.  For each move, there is a bijection of the colorings of $D$ and $D'$ such that the corresponding colorings are the identical outside $\Delta$.  When the move is of type R1, R2 or R3, it follows from that  $(X, R)$ is a biquandle \cite{FJK04}.  
When the move is of type V1, V2 or V3, it follows from that 
 $(X, V)$ is a biquandle and $V^2=1$. The case of V4 follows from the equality $(V \times 1)(1 \times V)(R \times 1)= (1 \times R)(V \times 1) (1 \times V)$.  
The case of type T1 follows from $V (T \times 1) = (1 \times T)V$. The case of type T2 follows from $T^2=1$.  The case of type T3 follows from $(T \times T)R (T \times R) = VRV$.
\qed

\vspace{0.3cm}
\indent
{\it Proof of Theorem~\ref{thm:bourgoinproductcoloring}}.  
Let ${\cal B}(Q, \ast) = (X, R, V, T)$.  
By Theorem~\ref{thm:coloring} it is sufficient to show that for a virtual link diagram $D$ there is a bijection between the set of $(X, R, V, T)$-colorings of $D$, denoted by ${\rm Col}(D, (X, R, V, T))$, and the Cartesian product of ${\rm Col^u}(D, (Q, \ast))$ and ${\rm Col^l}(D, (Q, \ast))$, where ${\rm Col^u}(D, (Q, \ast))$ is the set of upper $(Q, \ast)$-colorings of $D$ and ${\rm Col^l}(D, (Q, \ast))$ is the set of lower $(Q, \ast)$-colorings of $D$.  Since $D$ is a virtual link diagram, an  $(X, R, V, T)$-coloring of $D$  is nothing more than an $(X, R, V)$-coloring of $D$.  
Let $c: E(D) \to X$ be an $(X, R, V)$-coloring of $D$, where $E(D)$ is the set of edges of $D$.  Let $c^{\rm u}: E(D) \to Q$ and $c^{\rm l}: E(D) \to Q$ be the maps defined by $ c^{\rm u}= p_1 c$ and $c^{\rm l} = p_2 c$ where $p_i$ is the $i$th factor projection $X=Q^2 \to Q$.  Then $c^{\rm u}$ is  an upper $(Q, \ast)$-coloring of $D$, and $c^{\rm l}$ is a lower $(Q, \ast)$-coloring of $D$.  Conversely for any upper $(Q, \ast)$-coloring $c^{\rm u}$ and any lower $(Q, \ast)$-coloring $c^{\rm l}$, the map $c=(c^{\rm u}, c^{\rm l}): E(D) \to X$ is an $(X, R, V)$-coloring of $D$. \qed 

\vspace{0.3cm}
In order to prove Theorem~\ref{thm:twofoil}, we give a proposition on ${\cal B}(Q, \ast)$-colorings of the non-orientable virtual $m$-foil $F_m$.  

Let $(Q, \ast)$ be a quandle and ${\cal B}(Q, \ast)$ the standard twisted quandle.  
Consider a  ${\cal B}(Q, \ast)$-coloring of the diagram $F_m$  and let $x, y, z$ and $w$ be elements of $Q^2$ given the edges as depicted in Figure~\ref{fgtwistedmfoil}.  Put $x=(x_1, x_2)$ and $y=(y_1, y_2) \in Q^2$.  

\begin{proposition}\label{prop:coloringFm}
In this situation, $x_1$ and $y_1$ satisfy 
\begin{equation}
x_1 = x_1^{(y_1 x_1)^m} \quad \mbox{and} \quad 
y_1 = y_1^{(x_1 y_1)^m}. 
\end{equation}
And $x_2$ and $y_2$ are determined by 
\begin{equation} 
\left\{
\begin{array}{l}
x_2 = y_1^{(x_1 y_1)^n} \\ 
y_2 = x_1^{(y_1 x_1)^{n-1}y_1} 
\end{array}
\right. 
\quad  \mbox{if $m=2n$, and } 
\left\{
\begin{array}{l}
x_2 = x_1^{(y_1 x_1)^n y_1} \\ 
y_2 = y_1^{(x_1 y_1)^n} 
\end{array}
\right. 
\quad  \mbox{if $m=2n+1$.}
\end{equation}
Conversely, for any elements $x_1$ and $y_1$ of $Q$ satisfying $(1)$, there exists a unique ${\cal B}(Q, \ast)$-coloring of $F_m$. 
\end{proposition}

\vspace{0.3cm}
\indent
{\it Proof}.  
The standard twisted product of a quandle is given by 
\begin{equation*}
R( (x_1, x_2), (y_1, y_2)) = ((y_1, y_2^{x_2}), (x_1^{y_1}, x_2)).   
\end{equation*}
Applying that $m$ times, we have 
\begin{equation} 
(z, w) = (( x_1^{(y_1 x_1)^{n-1}y_1}, x_2^{(y_2 x_2)^n} ),  ( y_1^{(x_1 y_1)^n}, y_2^{ (x_2 y_2)^{n-1} x_2})) 
\quad  \mbox{if $m=2n$,} 
\end{equation}
and 
\begin{equation} 
(z, w) = (( y_1^{(x_1 y_1)^n}, y_2^{ (x_2 y_2)^n x_2}),  ( x_1^{(y_1 x_1)^n y_1}, x_2^{(y_2 x_2)^n} )) 
\quad  \mbox{if $m=2n+1$.}
\end{equation}
Since $(x, y)= V (T\times T)(z,w)$, we have 
\begin{equation} 
\left\{
\begin{array}{ll}
x_1 = y_2^{(x_2 y_2)^{n-1}x_2}, &  x_2 = y_1^{(x_1 y_1)^n} \\ 
y_1 = x_2^{(y_2 x_2)^n},             &   y_2 = x_1^{(y_1 x_1)^{n-1}y_1} 
\end{array}
\right. 
\quad  \mbox{if $m=2n$,} 
\end{equation}
and 
\begin{equation} 
\left\{
\begin{array}{ll}
x_1 = x_2^{(y_2 x_2)^n},             &  x_2 = x_1^{(y_1 x_1)^n y_1} \\ 
y_1 = y_2^{(x_2 y_2)^n x_2},      &   y_2 = y_1^{(x_1 y_1)^n} 
\end{array}
\right. 
\quad  \mbox{if $m=2n+1$.} 
\end{equation}
From (5) and (6) we have (2), and eliminating $x_2$ and $y_2$ we have (1).  Conversely, for any elements $x_1$ and $y_1$ of $Q$ satisfying $(1)$, let $x_2$ and $y_2$ be as in (2).  Then $x=(x_1, x_2)$ and $y=(y_1, y_2)$ satisfy (5) and (6).  Thus we have a coloring.  \qed 

\vspace{0.3cm}
For a quandle $(Q, \ast)$, let $\Delta_m(Q, \ast)$ denote a subset 
\begin{equation}
\{ (a,b) \in Q^2 \mid a= a^{(ba)^m}, \, b= b^{(ab)^m} \}
\end{equation}
of $Q^2$.  

\begin{corollary}\label{cor:colorFm}
For any quandle $(Q, \ast)$, assigning $(x_1, y_1)$ to a ${\cal B}(Q, \ast)$-coloring of $F_m$ as in Proposition~\ref{prop:coloringFm} 
is a bijection from ${\rm Col}(F_m, {\cal B}(Q, \ast))$, the set of  ${\cal B}(Q, \ast)$-colorings of $F_m$, to $\Delta_m(Q, \ast)$.  
\end{corollary}

\vspace{0.3cm}
\indent
{\it Proof of Theorem~\ref{thm:twofoil}}.  
(1) Let $(Q_{2m}, \ast)$ be the dihedral quandle of order $2m$, i.e., $Q_{2m}= \Z / 2m\Z$ and $a \ast b = 2b -a$.  
It is easily seen that $\Delta_m(Q_{2m}, \ast) = (Q_{2m})^2$.  Thus, by Corollary~\ref{cor:colorFm},  $\# {\rm Col}(F_m, {\cal B}(Q_{2m}, \ast)) = \# \Delta_m(Q_{2m}, \ast) = 4m^2$.  On the other hand, if $m > m' \geq 1$, then $\Delta_{m'}(Q_{2m}, \ast) \neq (Q_{2m})^2$, for $(1,0)\in (Q_{2m})^2$ does not belong to $ \Delta_{m'}(Q_{2m}, \ast)$.  
Hence, by Corollary~\ref{cor:colorFm} again, $\# {\rm Col}(F_{m'}, {\cal B}(Q_{2m}, \ast)) = \# \Delta_{m'}(Q_{2m}, \ast) < 4m^2$.  Therefore $F_m$ and $F_{m'}$ represent distinct twisted link. 

(2) Let $(Q_{2m+1}, \ast)$ be the dihedral quandle of order $2m+1$.  Since  $(1,0)\in (Q_{2m+1})^2$ does not belong to $ \Delta_{m}(Q_{2m+1}, \ast)$,  we have $\Delta_m(Q_{2m+1}, \ast) \neq (Q_{2m+1})^2$.   
By Corollary~\ref{cor:colorFm}, $\# {\rm Col}(F_{m}, {\cal B}(Q_{2m+1}, \ast)) = \# \Delta_{m}(Q_{2m+1}, \ast) < 
\# (Q_{2m+1})^2$.  
By Corollary~\ref{cor:bourgoinproductcoloring}, $F_m$ does not represent a virtual link.  \qed

\vspace{0.3cm}
{\it Acknowledgements.} 
The authors would like to express their gratitude to Roger Fenn for his valuable suggestions.   
The first author is partially supported by JSPS KAKENHI 22540093 and Grant-in-Aid for Research in Nagoya City University.  The second author  is partially supported by JSPS KAKENHI 21340015.

\end{document}